\documentclass{article}
\usepackage{amsmath}
\usepackage{amssymb}  % making IR, IN, etc.
\usepackage{epsfig}
\usepackage{amsthm}   % \begin{proof} \end{proof}
\usepackage[mathcal]{eucal}
\usepackage{tensor}
\usepackage{url} 

\newcommand{\F}{\mathcal{F}}

\newcommand{\R}{\mathbb{R}}
\newcommand{\BR}{\bar{\mathbb{R}}}

\DeclareMathOperator*{\cl}{cl}

\newcommand{\inner}[2]{\langle{#1},{#2}\rangle}

\newcommand{\pr}{\vdash}

\newcommand{\tos}{\rightrightarrows} % point-to-set mappings

\newtheorem{theorem}{Theorem}[section]
\newtheorem{lemma}[theorem]{Lemma}

\newtheorem{definition}{Definition}[section]

\title{Non-enlargeable operators and self-cancelling operators}

\author{B. F. Svaiter\thanks{ IMPA, Estrada Dona Castorina 110, 
    22460-320 Rio de
    Janeiro, Brazil ({\tt benar@impa.br}) }\hspace{.5em}
    \thanks{Partially supported by CNPq
    grants 300755/2005-8, 475647/2006-8 and by PRONEX-Optimization}
}

\begin{document}

\maketitle

\begin{abstract}
  The $\varepsilon$-enlargement of a maximal monotone operator is a
  construct similar to the Br{\o}ndsted and Rocakfellar
  $\varepsilon$-subdifferential enlargement of the subdifferential.
  Like the $\varepsilon$-subdifferential, the
  $\varepsilon$-enlargement of a maximal monotone operator has
  practical and theoretical applications.

  In a recent paper in Journal of Convex Analysis Burachik and Iusem
  studied conditions under which a maximal monotone operator is
  non-enlargeable, that is, its $\varepsilon$-enlargement coincides
  with the operator. Burachik and Iusem studied these non-enlargeable
  operators in reflexive Banach spaces, assuming the interior of the
  domain of the operator to be nonempty. In the present work, we remove
  the assumption on the domain of non-enlargeable operators and also
  present partial results for the non-reflexive case.
  \\
\end{abstract}

%\pagestyle{plain}

%%%%%%%%%%%%%%%%%%%%%%%%%%%%%%%%%%%%%%%%%%%%%%%%%%%%%%%%%%%%%%%%%%%%%%%%%

\section{Introduction}

Let $X$ be a real Banach space.  We use the notation $X^*$ for the
topological dual of $X$ and 
$\inner{\cdot}{\cdot}$ for the duality product
in $X\times X^*$:
\[
\inner{x}{x^*}=x^*(x).
\]
Whenever necessary, we will identify $X$ with its image under the
canonical injection of $X$ into $X^{**}$. 
A point-to-set operator $T:X\tos X^*$ is a relation on $X\times X^*$:
\[ 
 T\subset X\times X^* 
\]
and $x^*\in T(x)$ means $(x,x^*)\in T$. From now on, for $T:X\tos
X^*$, we define
\[ -T=\{(x,-x^*)\;|\; (x,x^*)\in T\}.
\]
so that $ -T:X\tos X^*$, $(-T)(x)=-(T(x))$.
An operator $T:X\tos X^*$ is {\it monotone} if
\[
\inner{x-y}{x^*-y^*}\geq 0,\forall (x,x^*),(y,y^*)\in T.
\]
and it is {\it maximal monotone}  if it
is monotone and maximal (with respect to the inclusion) in the family
of monotone operators of $X$ to $X^*$.
Maximal monotone operators in Banach spaces arise, for example, in
the study of PDE's, equilibrium problems and calculus of variations.

The \emph{conjugate} of $f:X\to\BR$ is $f^*:X^*\to \BR$,
\[
 f^*(x^*)=\sup_{x\in X} \inner{x}{x^*}-f(x).
\]
and the \emph{effective domain} of $f$ is
\[
 ed(f)=\{x\in X\;|\; f(x)<\infty\}.
\]
The subdifferential of $f$ is the point to set operator $\partial f:X\tos X^*$,
\[
 \partial f(x)=\{ x^*\;|\; f(y)\geq f(x)+\inner{y-x}{x^*},\;\forall y\}.
\]
In a paper
where many fundamental techniques were introduced, Rockafellar proved
that the subdifferential of a proper, convex, lower semicontinuous
function in a Banach space is maximal monotone~\cite{MR0262827}.
Rockafellar's proof relied on the
\emph{$\varepsilon$-subdifferential}, a concept introduced previously
by Br{\o}ndsted and Rockafellar~\cite{MR0178103}, which is defined as
follows, for $f:X\to\BR$:
\[
\partial_\varepsilon f(x)=\{ x^*\;|\; f(y)\geq
f(x)+\inner{y-x}{x^*}-\varepsilon,\;\forall y\}.
\] 
Note that $\partial f\subset \partial_\varepsilon f$, for any
$\varepsilon \geq 0$, and the inclusion may be proper if
$\varepsilon>0$.  Hence, the $\varepsilon$-subdifferential is an
``enlargement'' of the subdifferential.
It is easy to check that for any $\varepsilon>0$, the
$\varepsilon$-subdifferential of $f$ is non-empty at any point where
$f$ is finite. One of the key properties of the
$\varepsilon$-subdifferential used on Rockafellar proof is the fact,
proved by Br{\o}ndsted and Rockafellar~\cite{MR0178103}, that points at
$\partial _\varepsilon f$ are close to $\partial f$, and this distance
can be estimated. This property is know as Br{\o}ndsted-Rockafellar
property of the $\varepsilon$-subdifferential.
Although created by Br{\o}ndsted and Rockafellar for theoretical purposes,
the $\varepsilon$-subdifferential has extensive practical applications
in convex optimization~\cite{MR1110085,MR1453305,MR1247617,MR1354434,MR1295240}.

If $T:X\tos X^*$ is maximal monotone, then inclusion on $T$ may be
characterized by a family of inequalities:
\[
 (x,x^*)\in T\iff\Big( \inner{x-y}{x^*-y^*}\geq 0,\;
 \forall (y,y^*)\in T\Big).
\]
Martinez-Legaz and Thera~\cite{MR1422399} observed that the above
inequality could be relaxed, in order to define an enlargement of $T$.
Burachik, Iusem and
Svaiter proposed the $T^\varepsilon$ enlargement~\cite{MR1463929} as follows:
for $\varepsilon\geq 0$,
\begin{equation}
  \label{eq:teps}
  (x,x^*)\in T^\varepsilon\iff\Big( \inner{x-y}{x^*-y^*}\geq -\varepsilon,\;
 \forall (y,y^*)\in T\Big).
\end{equation}

The $T^\varepsilon$ enlargement has many similarities to the
$\varepsilon$-subdifferential proposed by Br{\o}ndsted and
Rockafellar~\cite{MR0178103}.  For example, in the interior of the
domain of $T$, for $\varepsilon$ bounded away from $0$, the mapping
\[
(x,\varepsilon)\mapsto T^\varepsilon(x)=\{ x^*\;|\;(x,x^*)\in T^\varepsilon\}
\]
is locally Lipschitz continuous, with respect to the Hausdorff 
metric.  This enlargement also satisfy (in reflexive spaces) a
property similar to the Br{\o}ndsted-Rockafellar property of the
$\varepsilon$-subdifferential.  Beside that, the $T^\varepsilon$
enlargement has also theoretical~\cite{MR1871464,MR1737322,MR2265756}
and algorithmic
applications~\cite{MR1756912,MR1682740,MR1737312,MR1783979,
MR1811750,MR1871872,MR1980655}. For
a survey in the subject, see \cite{MR2353163}.

Our aim is to investigate those maximal monotone operators $T:X\tos X^*$ which
are ``non-enlargeable'', that is, 
\begin{equation}
  \label{eq:theproblem}
  T^\varepsilon=T,\qquad \forall \varepsilon\geq 0.
\end{equation}
This question has been previously addressed by Burachik and
Iusem~\cite{MR2291554} and the present work is inspired in that article
of Burachik and Iusem.

It shall be noted that the $T^\varepsilon$ enlargement is one among a
family of enlargements, defined and studied on~\cite{MR1802238}.  These
enlargements share some basic properties and $T^\varepsilon$
is the biggest element in this family. Moreover, if $T$ happens to be
the subdifferential of some convex function $f$, then the
$\varepsilon$-subdifferential of $f$ also belongs to this family and
the inclusion
\[ 
\partial_\varepsilon f\subset (\partial f)^\varepsilon
\]
is proper, in general.

The $T^\varepsilon$ enlargement is closely tied to the Fitzpatrick
function, which we discuss next. To honor Fitzpatrick, we shall use $\varphi$,
the Greek ``f'', to denote Fitzpatrick function~\cite{MR1009594} associated with a maximal monotone operator $T:X\tos X^*$:
\begin{equation}
  \label{eq:f.fitz}
  \varphi_T(x,x^*)=\sup_{(y,y^*)\in T}
\inner{x}{y^*}+\inner{y}{x^*}-\inner{y}{y^*}\,.
\end{equation}
Observe that $ \varphi_T$ is convex, lower semicontinuous on the $w\times w^*$
topology of $X\times X^*$ and
\begin{equation}
  \label{eq:fitz.and.t}
 \varphi_T(x,x^*)\geq \inner{x}{x^*}, 
 \qquad T=\{ (x,x^*)\;|\; \varphi_T(x,x^*)=\inner{x}{x^*}\}  \,.
\end{equation}
The above inequality is a generalization of Fenchel-Young inequality. Indeed
if $f$ is a proper convex lower semicontinuous function on $X$, then
\[
 f(x)+f^*(x^*)\geq \inner{x}{x^*},\qquad
\partial f=\{ (x,x^*)\;|\; f(x)+f^*(x^*)=\inner{x}{x^*}\}.
\]
So defining $f$, defining 
the Fenchel-Young function associated with $f$,
\begin{equation}
  \label{eq:def.hfy}
h_{FY}:X\times X^*\to\BR,\qquad h_{FY}(x,x^*)=f(x)+f(x^*),  
\end{equation}
we have a convex function bounded bellow by the duality product and
equal to it at $\partial f$. Fitzpatrick proved that associated with
each maximal monotone operator  $T$ there is a family of functions with
these properties and that $\varphi_T$ is the minimal element of this family.
Br{\o}ndsted and Rockafellar observed that the $\varepsilon$-subdifferential can
be characterized by the function $h_{FY}$:
\begin{align}
\nonumber
  \partial _\varepsilon f&=\{ (x,x^*)\;|\; f(x)+f^*(x^*)\leq \inner{x}{x^*}
 +\varepsilon\}\\
  &=\{ (x,x^*)\;|\; h_{FY}(x,x^*)\leq \inner{x}{x^*}
 +\varepsilon\}.
\end{align}
Likewise, it is trivial to check that $\varphi_T$ characterizes the
$T^\varepsilon$ enlargement of a maximal monotone $T:X\tos X^*$:
\begin{equation}
  \label{eq:teps.2}
   T^\varepsilon =\{ (x,x^*)\;|\; \varphi_T(x,x^*)\leq
\inner{x}{x^*}+\varepsilon\}.
\end{equation}

Given a maximal monotone operator $T:X\tos X^*$,
Fitzpatrick defined~\cite{MR1009594} the family $\mathcal{F}_T$ as those
convex, lower semicontinuous functions in $X\times X^*$ which are
bounded below by the duality product and coincide with it at $T$:
\begin{equation}
  \label{eq:def.ft}
  \F_T=\left\{ h\in \BR^{X\times X^*}
    \left|
      \begin{array}{ll}
        h\mbox{ is convex and lower semicontinuous}\\
        \inner{x}{x^*}\leq h(x,x^*),\quad \forall (x,x^*)\in X\times X^*\\
        (x,x^*)\in T 
        \Rightarrow 
        h(x,x^*) = \inner{x}{x^*}
      \end{array}
    \right.
  \right\}.
\end{equation}
Fitzpatrick proved that $\varphi_T$ belongs to this family and it is
its minimal element.  Moreover, he also proved that if $h\in \F_T$
then $h$ represents $T$ in the following sense:
\[
 (x,x^*)\in T\iff h(x,x^*)=\inner{x}{x^*}.
\]
For the case of the subdifferential of a proper convex
lower semicontinuous function $f$, defining $h_{FY}$ as \eqref{eq:def.hfy}, 
\[ 
 h_{FY}\in \mathcal{F}_{\partial f}.
\]
Moreover, $h_{FY}$ is separable. It would be most desirable to find
separable elements in $\mathcal{F}_T$. Unfortunately, this family has
a separable element if and only if $T$ is a subdifferential~\cite{MR1974634}. Another
interesting property of $h_{FY}$ is that this function is a fixed
point of the mapping
\[ 
\mathcal{J}:{\BR}^{X\times X^*}\to{\BR}^{X\times X^*},
\qquad \mathcal{J}g(x,x^*)=g^*(x^*,x).
\]
Burachik and Svaiter observed the $\mathcal{F}_T$ is invariant under
$\mathcal{J}$~\cite{MR1934748} and Svaiter proved that there always exist a
fixed point of $\mathcal{J}$ in $\mathcal{F}_T$\cite{MR1999934}.
These fixed points has meet some applications in the study of PDE'S
under the attractive name
``self-dual''~\cite{MR2123037,MR2121900,MR2235885,MR2310692,MR2334832,MR2379461}
in the pioneering works of Ghoussoub.

\section{Non-enlargeable operators}
\label{sec:none}

Direct use of \eqref{eq:teps.2} shows that problem \eqref{eq:theproblem}
is equivalent to 
finding those maximal monotone operators
$T$ such that
\begin{equation}
  \label{eq:a}
  ed(\varphi_T)=T.
\end{equation}
It has been recently proved in~\cite{Alves2008} and in~\cite{Bauschke2008}, independently,
that if a maximal monotone is convex, then it is an affine subspace of
$X\times X^*$.  As $\varphi_T$ is convex, the above condition implies
that $T$ is convex.  Therefore, we can reduce our 
problem to finding
those   maximal monotone operators which are affine subspaces and
 satisfy \eqref{eq:theproblem}.

If $T\subset X\times X^*$ and $(x_0,x_0^*)\in X\times X^*$, defining
\begin{align*}
 T_0&=T-\{(x_0,x_0^*)\}\\
  &=\{(x-x_0,x^*-x_0^*)\;|\; (x,x^*)\in T\}.
\end{align*}
We have
\begin{equation}
  \label{eq:tr}
  (T_0)^\varepsilon=T^\varepsilon-\{(x_0,x_0^*)\}.
\end{equation}
So, we can restrict our attention tho those maximal monotone operators
which are subspaces of $X\times X^*$ and satisfy \eqref{eq:theproblem}, and the
general case will be obtained by translations of these subspaces.

Define, for $B\subset X\times X^*$ 
\begin{equation}
  \label{eq:def.ort}
 B^\pr=
 \{(y,y^*)\;|\; \inner{x}{y^*}+\inner{y}{x^*}=0,\qquad \forall
 (x,x^*)\in B\} \,. 
\end{equation}
Note that $B^\pr$ can be written in terms of the annihilator of a
family in $(X\times X^*)^*$:
\[ 
B^\pr=\mbox{}^a{ \{ (x^*,x)\;|\; (x,x^*)\in B\}}
\]
\begin{lemma}
  \label{lm:ort}
  If $T\subset X\times X^*$ is maximal monotone and a subspace, then
  \begin{enumerate}
  \item $ T^\pr\subset \{ (x,x^*)\;|\; \varphi_T(x,x^*)=0\}$,
  \item $ T\cap T^\pr=T\cap \{(x,x^*)\;|\; \inner{x}{x^*}=0\}$.
  \end{enumerate}
\end{lemma}
\begin{proof}
  To prove item 1, take $(x,x^*)\in T^\pr$. As $(0,0)\in T$, for any
  $(y,y^*)\in T$, $\inner{y}{y^*}\geq 0$. Therefore
   \begin{align*}
     \varphi_T(x,x^*)&=\sup_{(y,y^*)\in T} 
         \inner{x}{y^*}+\inner{y}{x^*}-\inner{y}{y^*}\\
      &=\sup_{(y,y^*)}-\inner{y}{y^*}=0
   \end{align*}

   To prove item 2, first use item 1 to obtain
   \[ 
   T\cap T^\pr\subset T\cap\{ (x,x^*)\;|\; \varphi_T(x,x^*)=0\}.
   \]  
   As $\varphi_T(x,x^*)=\inner{x}{x^*}$ for any $(x,x^*)\in T$, we
   conclude 
   \[ 
   T\cap T^\pr\subset T\cap \{(x,x^*)\;|\; \inner{x}{x^*}=0\}.
   \]
  To prove the other inclusion, take $(x,x^*)\in T\cap \{(x,x^*)\;|\;
  \inner{x}{x^*}=0\}$. As $T$ is a subspace, if $(y,y^*)\in T$, then, for
  any $\lambda\in \R$
  \[
  \lambda(x,x^*)+(y,y^*)\in T.
  \]
  As $(0,0)\in T$,
  \[
  \inner{\lambda x+y}{\lambda x^*+y^*}=\lambda [
  \inner{x}{y^*}+\inner{y}{x^*}]
  +\inner{y}{y}\geq 0.
  \]
  As $\lambda$ is arbitrary, we conclude that the expression inside
  the brackets must be $0$.
\end{proof}
It is interesting to observe that $\inner{x}{x^*}$ is  non-linear
and non-convex in $(x,x^*)$. Even though, the points at $T$ where this
expression vanish is a subspace, which may be empty also.

We will be concerned with a special type of linear point-to set
operators
\begin{definition}
  \label{def:anti}
  An operator $A:X\tos X^*$ is \textit{\bf self-cancelling} if $A$ is
  a subspace and
  \[
  \inner{x}{x^*}=0,\qquad \forall (x,x^*)\in A.
  \]
\end{definition}
This definition is an extension of the definition of skew-symmetric
operators of Burachik and Iusem~\cite{MR2291554} and of the definition
of skew linear of Bauschke, Wang and Yau~\cite{Bauschke2008}. The
relations between these classes will be discussed in the
Section~\ref{sec:as}.

\begin{lemma}
  \label{lm:tech}
  If $A\subset X\times X^*$ is self-cancelling, then
  $
  A\subset A^\pr
  $.
\end{lemma}
\begin{proof}
  Take $(x,x^*),(y,y^*)\in A$. Then, $(x+y,x^*+y^*)\in A$ and so
  \[ 
  \inner{x+y}{x^*+y^*}=\inner{x}{x^*}+\inner{x}{y^*}+\inner{y}{x^*}+
  \inner{y}{y^*}=0.
  \]
  To end the proof, note that  $\inner{x}{x^*}=\inner{y}{y^*}=0$. 
\end{proof}
\begin{lemma}
  \label{lm:tech.2}
  If $A\subset X\times X^*$ is self-cancelling and $A^\pr$ is
  maximal monotone then, for any $(x_0,x_0^*)\in X\times X^*$, 
 the operator
  \[ 
  T=A^\pr+\{(x_0,x_0^*)\}
  \]
  is non-enlargeable, or equivalently, $ed(\varphi_T)=T$.
\end{lemma}
\begin{proof}
  In view of \eqref{eq:tr}, it suffices to prove this lemma for $(x_0,x_0^*)=0$.
  In that case, if $(x,x^*)\notin A^\pr$, there exists
  $(y,y^*)\in A$ such that
  \[
  \inner{x}{y^*}+\inner{y}{x^*}\neq 0.
  \]
  As $A$ is a subspace and $A\subset A^\pr$,
  \[
  \varphi_{A^\pr}(x,x^*)\geq 
  \sup_{\lambda \in \R }  \inner{x}{\lambda y^*}+\inner{\lambda y}{x^*}
  -\inner{\lambda y}{\lambda y^*}
  = \sup_{\lambda \in \R }  \inner{x}{\lambda y^*}+\inner{\lambda y}{x^*}.
  \]
  Combining the above equation we obtain $ \varphi_{A^\pr}(x,x^*)\geq\infty$.
\end{proof}
\begin{theorem}
  \label{th:main}
  If $X$ is reflexive then $T$ maximal monotone is non-enlargeable, if
  and only if there exists an self-canceling $A$ and $(x,x^*)\in
  X\times X^*$ such that $A^\pr$ is maximal monotone and
  \[ 
  T=A^\pr+\{(x,x^*)\}
  \]
  Moreover, if $T$ is non-enlargeable and $(x,x^*)\in T$, then the
   maximal $A$ satisfying the above condition is 
   \[
   (T-{(x,x^*)})^\pr.
   \]
\end{theorem}
\begin{proof}
  First assume that $T$ is non-enlargeable and $(0,0)\in T$.  Define
  \[ 
  A=T^\pr
  \]
  Using Lemma~\ref{lm:ort}, item 1, we conclude that $\varphi_T(x,x^*)=0$ for
  all $(x,x^*)\in A$. As $ed(\varphi_T)=T$, we conclude that $A\subset
  T$. Therefore,
  \[
  A=T\cap A.
  \]
  Combining the above equation with the definition of
  $A$ and Lemma~\ref{lm:ort}, item 2, we conclude that $A$ is
  self-cancelling. Moreover $A$ is the maximal self-cancelling operator
  contained in $T$.  As $T$ is a closed subspace and $X$ is reflexive,
  direct use of Hahn-Banach yields
  \[
  T=(T^\pr)^\pr=A^\pr.
  \]
  Note also that the above defined $A$ is maximal in
  the family
  \[
  \{ B\subset X\times X^*\;|\; T=B^\pr\}
  \]
  Conversely, if for some self-cancelling $A$, $T=A^\pr$, then according to 
  Lemma~\ref{lm:tech.2} $T$ is non-enlargeable.

  The general case follows now using~\eqref{eq:tr}. 
\end{proof}

\section{On maximal monotone operators obtained self-cancelling operators}
\label{sec:as}

Now we shall analyze those maximal monotone operators discussed in
Theorem~\ref{th:main}.
\begin{lemma}
  \label{lm:as.1}
  If $A\subset X\times X^*$ is self-cancelling and $A^\pr$ is monotone,
  then  $A^\pr$ is maximal monotone.
\end{lemma}
\begin{proof}
  Take $(x_0,x_0^*)\notin A^\pr$. Then there exists $(y,y^*)\in A$
  such that
  \[
  \inner{x_0}{y}+\inner{y}{x_0^*}\neq 0.
  \]
  Then, for any $\lambda\in \R$,
  \[
  \inner{x_0+\lambda y}{x_0^*+\lambda y^*}=\inner{x_0}{x_0^*}
  +\lambda[ \inner{x_0}{y^*}+\inner{y}{x_0^*}].
  \]   
  Combining the two above equations with the fact that $A$ is a
  subspace, we conclude that $\{(x_0,x_0^*)\}\cup A$ is not monotone. 
\end{proof}
Observe that in Lemma~\ref{lm:ort}, the  maximal monotone operator $T$
may not be on the family of $A^\pr$ with $A$ self-cancelling.
\begin{lemma}
  \label{lm:as.2}
  If $A$ is self-cancelling and $A^\pr$ is monotone, then 
  \[
   \tilde A= \cl_{w\times w*}A=\cl_{s\times w*} A
  \]
  is a maximal element of the family of self-cancelling operators, where
  $\cl_{w\times w*}$ and $\cl_{s\times w*} $ denotes the closure in the
  weak$\times$weak-$*$ and  strong$\times$weak-$*$ topologies, respectively.
\end{lemma}
\begin{proof}
  Suppose that $B =[\tilde A\cup\{x_0,x_0^*\}]$ is antisymmetric. In that case,
  \[
  B^\pr \subset A^\pr=(\tilde A)^\pr.
  \]
  In particular, $ B^\pr$ is monotone. Hence $B^\pr$ is maximal
  monotone and the above inclusion holds as an equality. As $\tilde A$
  is $w\times w*$ closed, $(x_0,x_0^*)\in \Tilde A$.
\end{proof}

A natural question is  whether  $A^\pr$ is maximal monotone
whenever $A$ is maximal self-cancelling.
Up to now we have a partial answer to this question.
\begin{lemma}
  \label{lm:sc.m}
  If $A:X\tos X^*$ is maximal self-cancelling, the $A^\pr$ or
  $-A^\pr$ is maximal monotone.
\end{lemma}
\begin{proof}
  Recall that $ -A=\{(x,-x^*)\;|\; (x,x^*)\in A\}$ so that
 \[
  (-A)^\pr=-(A^\pr)\;,
 \]
  Take $(x,x^*),(y,y^*)\in A^\pr$. Suppose that
  \begin{equation}
    \label{eq:test}
      \inner{x}{x^*}>0,\qquad \inner{y}{y^*}<0.
  \end{equation}
  Then, for some $\theta\in(0,1)$
  \[
  z_\theta=\theta x+(1-\theta)y,\qquad 
  z_\theta^*=\theta x^*+(1-\theta)y^*.
  \]
  satisfy 
  \[
   \inner{z_{\theta}}{z^*_{\theta}}=0
  \]
  Hence, $(z_\theta,z^*_\theta)\in A$ and
  $\inner{z_\theta}{y^*}+\inner{y}{z_\theta^*}=0$.  Direct use of the
  definition of $z_\theta,z_\theta^*$ gives
  \[ 
  \theta(x,x^*)= (z_\theta,z^*_\theta)-(1-\theta)(y,y^*),
  \]
  which readily implies
  \begin{align*}
    \theta^2\inner{x}{x^*}&=
     \inner{z_\theta-(1-\theta) y}{z_\theta^*-(1-\theta) y^*} \\
    &=(1-\theta)^2\inner{y}{y^*}
  \end{align*}
  in contradiction with~\eqref{eq:test}. Therefore \eqref{eq:test} can
  not hold for $(x,x^*),(y,y^*)\in A^\pr$ and $ A^\pr$ or
  $-A^\pr$ is monotone. Maximal monotonicity of $A^\pr$ or $-A^\pr$
  now follows from Lemma~\ref{lm:as.1}
\end{proof}

Working in the setting of reflexive Banach spaces, Burachik and
Iusem~\cite{MR2291554}
defined a skew-symmetric operator as a linear continuous operator
$L:X\to X^*$ such that
\[ 
L=-L^*
\]
where $L^*$ is the adjoint of $L$. As $L^*:X^{**}\to X^*$, it is
natural to consider, in a reflexive Banach space, $L^*:X\to X^*$. In that
case,
$L^*$ is defined as
\[
 \inner{Lx}{y}=\inner{x}{L^*y},\qquad \forall x,y\in X.
\]
Note that $L^*=(-L)^\pr$.  Bauschke, Wang and
Yao~\cite{Bauschke2008}, still working in a reflexive Banach spaces,
extended this definition of adjoint to an arbitrary linear
point-to-set operator $L:X\tos X^*$ as $L^*=(-L)^\pr$. For these
authors, a point to set operator $L:X\tos X^*$ is skew if it is linear
and $L=-L^*$. It is trivial to verify that a skew-symmetric
operator (in the sense of \cite{MR2291554}) is always a skew operator
(in the sense on \cite{Bauschke2008}).
\begin{lemma}
  Let $M:X\tos X^*$
  be a linear point to set operator.
  \begin{enumerate}
  \item If $M$ is a skew operator, then it is maximal self-cancelling.
  \item If $M$ is maximal self-cancelling and $D(M)$  is
    closed  then it is skew.
  \item If $M$ is maximal self-cancelling,  $R(M)$ is
    closed and $X$ is reflexive, then it is skew.
  \end{enumerate}
\end{lemma}
\begin{proof}
  To prove item 1, suppose that $M$ is skew. Then, $M$ is
  self-cancelling. If $A$ is self-cancelling and $M\subset A$, then
  using Lemma~\ref{lm:tech} and \eqref{eq:def.ort} we have
  \[
  A\subset A^\pr\subset M.
  \]
  Therefore, $A=M$.

  To prove item 2, suppose that $M$ is maximal self-cancelling. Take
  \[
  (x_0,x_0^*)\in M^\pr
  \]
  If $x_0\notin D(M)$, then there exists $y^*$ such that
  \[
  \inner{x}{y^*}=0,\;\;\forall x \in D(M),\qquad
  \inner{x_0}{y^*}\neq 0.
  \]
  In that case, $(0,y^*)\in M$ and $ \inner{x_0}{y^*}+\inner{0}{x_0^*}\neq 0$,
  in contradiction with the assumption $  (x_0,x_0^*)\in M^\pr$.
  Hence, $x_0\in D(A)$ and there exists $z^*$ such that $(x_0,z^*)\in M$.
  Therefore,
  \[
  (x_0,x_0^*)-(x_0,z^*)\in M^\pr.
  \]
  To simplify the notation, let $u^*=x_0^*-z^*$. We have just proved that
  $(0,u^*)\in M^\pr$. Let
  \[
  V=\mathrm{span}(M\cup\{(0,u^*)\}).
   \]
  If $(x,x^*)\in M$ and $\lambda\in \R$, then
  \begin{align*}
    \inner{x}{x^*+\lambda u^*}&=\inner{x}{x^*}+\lambda\inner{x}{u^*}=0.
  \end{align*}
  Hence, $V$ is self-cancelling. As $M$ is maximal self-cancelling,
  \[
  (0,u^*)=(x_0,x_0^*)-(x_0,z^*)\in M,
  \]
  and $(x_0,x_0^*)\in M$. Altogether, using also Lemma~\ref{lm:tech} we have
  \[
  M\subset M^\pr\subset M
  \]
  and so $M$ is skew.

  Item 3 follows  from item 2, applied to $X'=X^*$ and 
  \[
  M'=\{ (x^*,x)\;|\; (x,x^*)\in M\}.
  \]
\end{proof}

\section{Acknowledgments}

The author thanks M. Marques Alves for the reading of the first draft
of this work and for the positive criticism and advises.

\bibliographystyle{plain} 

%\bibliography{monotone,prox,epsdiff}

%\def\cprime{$'$}

\end{document}